\documentclass[11pt]{amsart}

\usepackage{amsmath, amsthm}
\usepackage{amssymb}
\usepackage{amsfonts}
\usepackage{latexsym}
\usepackage{mathrsfs}
\usepackage{bm}
\usepackage{graphicx}
\usepackage[dvipsnames,usenames]{color}

\newcommand{\be}{\begin{equation}}
\newcommand{\ene}{\end{equation}}
\newcommand{\br}{\begin{remark}}
\newcommand{\er}{\end{remark}}
\newcommand{\bl}{\begin{lem}}
\newcommand{\el}{\end{lem}}
\newcommand{\bcor}{\begin{cor}}
\newcommand{\ecor}{\end{cor}}
\newcommand{\bpro}{\begin{pro}}
\newcommand{\epro}{\end{pro}}
\newcommand{\ben}{\begin{enumerate}}
\newcommand{\een}{\end{enumerate}}
\newcommand{\bp}{\begin{proof}}
\newcommand{\ep}{\end{proof}}
\newcommand{\bpo}{\begin{pro}}
\newcommand{\epo}{\end{pro}}
\newcommand{\beq}{\begin{equation*}}
\newcommand{\eeq}{\end{equation*}}
\newcommand{\bear}{\begin{eqnarray}}
\newcommand{\eear}{\end{eqnarray}}
\newcommand{\beqar}{\begin{eqnarray*}}
\newcommand{\eeqar}{\end{eqnarray*}}
\newcommand{\brem}{\begin{remark*}}
\newcommand{\erem}{\end{remark*}}
\newcommand{\bt}{\begin{theorem}}
\newcommand{\et}{\end{theorem}}
\newcommand{\bDe}{\begin{Def}}
\newcommand{\eDe}{\end{Def}}

\newcommand{\Z}{\mathbb{Z}}
\newcommand{\R}{\mathbb{R}}
\newcommand{\C}{\mathbb{C}}
\newcommand{\Cc}{\mathcal{C}_\gamma}

\DeclareMathOperator{\re}{Re}
\renewcommand{\Re}{\re}

\newcommand{\sech}{\operatorname{sech}}

\renewcommand{\tilde}[1]{\widetilde{#1}}

\newcommand{\eps}{\varepsilon}

\newcommand{\wep}{Weil-Petersson }

\newcommand{\sT}{\mathcal{T}}

\newcommand{\sM}{\mathcal{M}}
\newcommand{\lsys}{\ell_{sys}(X)}
\newcommand{\eo}{\epsilon_2}
\newcommand{\eoh}{\overline{\epsilon}_2}

\renewcommand{\leq}{\leqslant}
\renewcommand{\geq}{\geqslant}

\DeclareMathOperator{\Vol}{Vol_{WP}}

\DeclareMathOperator{\inj}{inj}

\DeclareMathOperator{\arcsinh}{arcsinh}

\DeclareMathOperator{\Teich}{Teich}

\DeclareMathOperator{\Diff}{Diff}

\newcommand{\Mod}{\mbox{\rm Mod}}
\DeclareMathOperator{\Ric}{Ric^{WP}}
\DeclareMathOperator{\HolK}{HolK^{WP}}
\DeclareMathOperator{\Sca}{Sca^{WP}}

\newcommand{\param}{{\mathchoice{\mkern1mu\mbox{\raise2.2pt\hbox{$\centerdot$}}\mkern1mu}{\mkern1mu\mbox{\raise2.2pt\hbox{$\centerdot$}}\mkern1mu}{\mkern1.5mu\centerdot\mkern1.5mu}{\mkern1.5mu\centerdot\mkern1.5mu}}}

\numberwithin{equation}{section}

\theoremstyle{plain}
\newtheorem{theorem}{Theorem}[section]

\newtheorem{corollary}[theorem]{Corollary}
\newtheorem{lemma}[theorem]{Lemma}
\newtheorem{proposition}[theorem]{Proposition}

\newtheorem*{rem*}{Remark}

\theoremstyle{definition}

\newtheorem{remark}[theorem]{Remark}
\theoremstyle{definition}
\newtheorem*{remarksenv}{Remarks}

\begin{document}

\title[Uniform Bounds]
{Uniform bounds on harmonic Beltrami differentials and Weil-Petersson curvatures}

\author{Martin Bridgeman}
\address[M. ~B. ]{Boston College, Chestnut Hill, Ma 02467, USA}
\email{bridgem@bc.edu}

\author{Yunhui Wu}
\address[Y. ~W. ]{Tsinghua University, Haidian District, Beijing 100084, China}
\email{yunhui\_wu@mail.tsinghua.edu.cn}

\begin{abstract}
In this article we show that for every finite area hyperbolic surface $X$ of type $(g,n)$ and any harmonic Beltrami differential $\mu$ on $X$, then the magnitude of  $\mu$ at any point of small injectivity radius is uniform bounded from above by the ratio of the \wep norm of $\mu$ over the square root of the systole of $X$ up to a uniform positive constant multiplication.

We  apply the uniform bound above to show that the \wep Ricci curvature, restricted at any hyperbolic surface of short systole in the moduli space, is uniformly bounded from below by the negative reciprocal of the systole up to a uniform positive constant multiplication. As an application, we show that the average total \wep scalar curvature over the moduli space is uniformly comparable to $-g$ as the genus $g$ goes to infinity. 
\end{abstract}

\subjclass{30F60, 53C21, 32G15}
\keywords{Uniform bounds, harmonic Beltrami differentials, Weil-Petersson curvature, total scalar curvature}

\maketitle

\section{Introduction}
In this paper, we derive uniform bounds on the curvature of the Weil-Petersson metric on $\sM^n_g$ the moduli space of conformal structures on the  surface of genus $g$ with $n$ punctures where $3g+n\geq 5$. We write $\sM_g$ for $\sM^0_g$ for simplicity. These bounds depend on new uniform bounds for the norm of harmonic Beltrami differentials in terms of injectivity radius. 

Let $X\in \sM^n_g$. Recall that the \emph{systole} $\lsys$ of $X$ is shortest length of closed geodesics in the hyperbolic surface $X$ and for $z \in X$, the {\em injectivity radius} $\inj(z)$ is the maximum radius of an embedded ball centered at $z$. We denote the Margulis constant in dimension two by 
$$\eo= \sinh^{-1}(1).$$ By the Collar Lemma, for $r(z) \leq \eo$, then $z$ is either contained in a  collar $\mathcal C_\gamma$ about a closed geodesic $\gamma$ or $z$ is in a  neighborhood $\mathcal C_c$ about a cusp $c$. The tangent space $T_X\sM^n_g$ of $\sM^n_g$ at $X$ can be identified with the space of harmonic Beltrami differentials on $X$. Let $\mu \in T_X\sM^n_g$. We denote by $||\mu||_{WP}$ the \wep norm of $\mu$, which is also the $L^2$-norm of $\mu$ on $X$. One consequence of our analysis is the following Proposition.

\begin{proposition}  \label{lnubw-0}
 Let $X\in \sM^n_g$ with $\lsys \leq 2\eo$.
 Then for any $\mu \in T_X\sM^n_g$  a harmonic Beltrami differential and $z \in X$  with injectivity radius $\inj(z)\leq \eo$ 
\[|\mu(z)|^2\leq   \frac{||\mu ||_{WP}^2}{\inj(z)} \leq 2\frac{||\mu||_{WP}^2}{ \ell_{sys}(X)}.\]
\end{proposition}

\br
In \cite[Corollary 11]{Wolpert5}, Wolpert proved a similar bound when $\lsys$ is smaller than a positive constant depending on $g$ and $n$. Our approach is similar to Wolpert's, but using a detailed analysis of the thin parts, we are able to obtain the above uniform bounds independent of $g$ and $n$. Actually we will prove certain more precise uniform bounds which are Proposition \ref{fourier1} and Lemma \ref{up-cusp}.  One may see Section \ref{ub-hbd} for more details.
\er

Using Proposition \ref{lnubw-0}, we derive uniform lower bounds on   \wep curvatures. More precisely, we prove
\bt \label{i-LUB-Ric}
For any $X \in \sM^n_g$ with $\lsys \leq2\eo $, then
\ben
\item for any $\mu \in T_X\sM^n_g$ with $||\mu||_{WP}=1$, the \wep Ricci curvature satisfies that
\[\Ric(\mu)\geq -\frac{4}{\lsys}.\] 

\item  The \wep scalar curvature at $X$ satisfies that
\[\Sca(X) \geq -\frac{4}{\lsys} \cdot (3g-3+n).\]
\een 
\et

\br
In \cite{Teo09} Teo showed that for any $X \in \sM_g$, 
\ben
\item $\Ric \geq -2 C(\frac{\lsys}{2})^2.$

\item $\Sca(X) \geq -(6g-6)C(\frac{\lsys}{2})^2.$
\een
Here the function $C(\cdot)$ is given by \eqref{defn:C-inj}. As the systole $\ell_{sys}(X)$ of $X$ tends to zero,  $C(\frac{\lsys}{2})^2 = \frac{4}{\pi\ell_{sys}(X)^2} + o(\frac{1}{\ell_{sys}(X)^2})$. Also $C(\frac{\lsys}{2})^2$ tends to $\frac{3}{4\pi}$ as $\ell_{sys}(X)$ goes to infinity. Compared to Teo's result, we obtain a better growth rate as $\lsys \to 0$. Actually this growth rate $\frac{-1}{\lsys}$ is optimal: Wolpert in \cite[Theorem 15]{Wolpert5} or \cite[Corollary 16]{Wolpert5} computed the \wep holomorphic sectional curvature along the gradient of certain geodesic length function and showed that it behaves as $\frac{-3}{\pi \ell_{\alpha}}+O(\ell_\alpha)$ as $\ell_\alpha \to 0$, where $\alpha \subset X$ is a nontrivial loop. Part (1) of Teo's results above in particular implies that the \wep sectional curvature, restricted on any $\epsilon$-thick part of the moduli space, is uniformly bounded from below by a negative constant only depending on $\eps$. This was first obtained by Huang in \cite{Huang}. One may also see \cite{Wolf-W-1} for more general statements. 
\er

\br
The assumption $\lsys \leq2\eo$ in Theorem \ref{i-LUB-Ric} can \emph{not} be removed. One may see this in the following two different ways; (1). Tromba \cite{Tromba86} and Wolpert \cite{Wol86} showed that for all $X\in \sM_g$, $$\Sca(X)\leq \frac{-3}{4\pi}\cdot (3g-2).$$ In particular for large enough $g$, the uniform lower bound for scalar curvature in Theorem \ref{i-LUB-Ric} does not hold for Buser-Sarnak surface $\mathcal{X}_g$ (see \cite{BS94}) whose injectivity radius grows like $\ln{(g)}$ as $g\to \infty$. Similarly for (2). It was shown in \cite[Theorem 1.1]{Wolf-W-1} that if $\lsys$ is large enough, then $$\min_{\emph{span}\{\mu,v\} \subset T_X\sM_g} K^{WP}(\mu,v)\leq -C<0$$ where $C>0$ is a uniform constant independent of $g$. In particular, the uniform lower bound for Ricci curvature in Theorem \ref{i-LUB-Ric} does not hold for Buser-Sarnak surface $\mathcal{X}_g$ in \cite{BS94} for large enough $g$.
\er

Let $X \in \sM^n_g$ with $\lsys \leq 2\eo$, and let $P(X)\subset T_X \sM^n_g$ be the linear subspace generated by the gradient of short closed geodesic length functions and $P(X)^\perp$ be its perpendicular. One may see \eqref{def-P} and \eqref{def-P-p} for the precise definitions. Our next result says that the \wep curvature along any plane in $T_X\sM_g$ containing a $\mu \in P(X)^\perp$ is uniformly bounded from below. More precisely,  
\bt \label{i-ub-orth}
Let $X \in \sM^n_g$ with $\lsys \leq 2\eo$, then for any $\mu\neq 0 \in P(X)^\perp$ and $v\in T_X \sM^n_g$, the \wep sectional curvature $K^{WP}(\mu,v)$ along the plane spanned by $\mu$ and $v$ satisfies that
\beqar \label{ulb-K}
K^{WP}(\mu,v)\geq -4.
\eeqar
\et
\noindent It would be \emph{interesting }to find upper bounds for $K^{WP}(\mu,v)$ in terms of certain measurements of $\mu$ and $v$.\\

Recall that the boundary $\partial \sM_g$ of $\sM_g$ consists of nodal surfaces. As $X$ goes to $\partial \sM_g$, the \wep scalar curvature $\Sca(X)$ always blows up to $-\infty$ because the \wep sectional curvature at $X$ along certain direction goes to $-\infty$ (e.g., see \cite{Schumacher} or \cite[Corollary 16]{Wolpert5}). It was not known  whether the total scalar curvature $\int_{\sM_g}\Sca(X) dX$ is finite. We will show it is truly finite. Moreover, combining Theorem \ref{i-LUB-Ric} and a result of Mirzakhani in \cite{MM4} we will determine the asymptotic behavior of $\int_{\sM_g}\Sca(X) dX$ as $g\to \infty$. More precisely, we prove
\bt \label{i-total}
As $g\to \infty$,
\[\frac{\int_{\sM_g} \Sca(X)dX}{\Vol(\sM_g)} \asymp -g.\]
\et
\

\noindent \textbf{Notation.} In this paper, we say two functions $$f_1(g) \asymp f_2(g)$$ if there exists a universal constant $C\geq 1$, independent of $g$, such that $$\frac{f_2(g)}{C} \leq f_1(g) \leq C f_2(g).$$\

\noindent \textbf{Plan of the paper.} Section \ref{s-np} provides some necessary background and the basic properties on Teichm\"uller theory and the Weil-Petersson metric. Refined results of Proposition \ref{lnubw-0} are proved in Section \ref{ub-hbd}. We prove several results on uniform lower bounds for \wep curvatures including Theorem \ref{i-LUB-Ric} and \ref{i-ub-orth}. Theorem \ref{i-total} is proved in Section \ref{s-WP-scal}.\\

\noindent \textbf{Acknowledgements.}
The authors would like to thank Jeffrey Brock, Ken Bromberg and Michael Wolf for helpful conversations on this project.

\section{Preliminaries}\label{s-np}
In this section, we set our notation and  review the relevant background material on Teichm\"uller space and \wep curvature. 
\subsection{Teichm\"uller space.} \label{sec:wp background}
We denote by $S^n_{g}$ an oriented surface of genus $g$ with $n$ punctures where $3g+n\geq 5$. Then the Uniformization theorem implies that the surface $S^n_g$ admits hyperbolic metrics of constant curvature $-1$. We let $\sT^n_g$ be the Teichm\"uller space of surfaces of genus $g$ with $n$ punctures, which we consider as the equivalence classes under the action of the group $\Diff_0(S^n_g)$ of diffeomorphisms isotopic to the identity of the space of hyperbolic surfaces $X=(S^n_g,\sigma(z)|dz|^2)$. The tangent space $T_X\sT^n_g$ at a point $X=(S^n_g,\sigma(z)|dz|^2)$ is identified with the space of finite area {\it harmonic Beltrami differentials} on $X$, i.e. forms on $X$ expressible as 
$\mu=\overline{\psi}/\sigma$ where $\psi \in Q(X)$ is a holomorphic quadratic differential on $X$. Let $z=x+iy$ and $dA=\sigma(z)dxdy$ be the volume form. The \textit{Weil-Petersson metric} is the Hermitian
metric on $\sT_g$ arising from the the \textit{Petersson scalar  product}
\begin{equation}
 \left<\varphi,\psi \right>= \int_X \frac{\varphi \cdot \overline{\psi}}{\sigma^2}dA\nonumber
\end{equation}
via duality. We will concern ourselves primarily with its Riemannian part $g_{WP}$. Throughout this paper we denote by $\Teich(S^n_g)$ the Teichm\"uller space endowed with the Weil-Petersson metric. By definition it is easy to see that the mapping class group $\Mod_g^n:=\Diff^+(S^n_g)/\Diff^0(S^n_g)$ acts on $\Teich(S^n_g)$ as isometries. Thus, the \wep metric descends to a metric, also called the \wep metric, on the moduli space of Riemann surfaces $\sM_g^n$ which is defined as $\sT^n_g/\Mod_g^n$. Throughout this paper we also denote by $\sM^n_g$ the moduli space endowed with the Weil-Petersson metric and write $\sM_g = \sM^0_g$ for simplicity.  One may refer to \cite{Wolpert4} for recent developments on \wep geometry.

\subsection{Weil-Petersson curvatures.} The \wep metric is K\"ahler. The curvature tensor of the \wep metric is given as follows. Let $\mu_{i},\mu_{j}$ be two elements in the tangent space $T_X\sM^n_g$ at $X$, so that the metric tensor  written in local coordinates is
\begin{eqnarray*}
g_{i \overline{j}}=\int_X \mu_{i} \cdot  \overline{\mu_j} dA.
\end{eqnarray*} 

\noindent For the inverse of $(g_{i\overline{j}})$, we use the convention
\begin{eqnarray*}
g^{i\overline{j}} g_{k\overline{j}}=\delta_{ik}.
\end{eqnarray*}

\noindent Then the curvature tensor is given by
\begin{eqnarray*}
R_{i\overline{j}k\overline{l}}=\frac{\partial^2}{\partial t^{k}\partial \overline{t^{l}}}g_{i\overline{j}}-g^{s\overline{t}}\frac{\partial}{\partial t^{k}}g_{i\overline{t}}\frac{\partial}{\partial \overline{t^{l}}}g_{s\overline{j}}.
\end{eqnarray*}

We now describe the curvature formula  of Tromba \cite{Tromba86} and Wolpert \cite{Wol86} which gives the curvature in terms of the Beltrami-Laplace operator $\Delta$. It has been applied to study various curvature properties of the Weil-Petersson metric. Tromba \cite{Tromba86} and Wolpert \cite{Wol86} showed that $\sM_g^n$ has negative sectional curvature. In \cite{Schumacher} Schumacher showed that $\sM_g^n$ has strongly negative curvature in the sense of Siu. Liu-Sun-Yau in \cite{LSY} showed that $\sM_g^n$ has dual Nakano negative curvature, which says that the complex curvature operator on the dual tangent bundle is positive in some sense. The third named author in \cite{Wu14-co} showed that the $\sM_g^n$ has non-positive definite Riemannian curvature operator. One can also see \cite{Huang-ajm, Huang, Teo09, Wol08, Wolpert5, Wu-IMRN} for other aspects of the curvature of $\sM_g^n$.

Set $D=-2(\Delta-2)^{-1}$ where $\Delta$ is the Beltrami-Laplace operator on $X=(S,\sigma|dz|^2) \in \sM_g^n$. The operator $D$ is positive and self-adjoint. 

\begin{theorem}[Tromba \cite{Tromba86}, Wolpert \cite{Wol86}]\label{cfow} 
The curvature tensor satisfies

\[R_{i\overline{j}k\overline{l}}=\int_{X} D(\mu_{i}\mu_{\overline{j}})\cdot (\mu_{k}\mu_{\overline{l}}) dA+\int_{X} D(\mu_{i}\mu_{\overline{l}})\cdot (\mu_{k}\mu_{\overline{j}}) dA.\]

\end{theorem}

\subsubsection{Weil-Petersson holomorphic sectional curvatures.}\label{subsec:hol-sec-curv}
Recall that a holomorphic sectional curvature is a sectional curvature along a holomorphic line. Let $\mu \in T_X\sM_g^n$ be a harmonic Beltrami differential. By Theorem \ref{cfow} the holomorphic sectional curvature $\HolK(\mu)$ along the holomorphic line spanned by $\mu$ is
\[\HolK(\mu)=\frac{-2\cdot \int_{X} D(|\mu|^2)\cdot (|\mu|^2) dA}{||\mu||_{WP}^4}.\]
Assume that $||\mu||_{WP}=1$. From \cite[Proposition 2.7]{Wolf-W-1}, which relies on an estimation of Wolf in \cite{Wolf12}, we know that 
\begin{equation}\label{HolK-eq}
-2 \int_{X} |\mu|^4 dA \leq \HolK(\mu)\leq -\frac{2}{3} \int_{X} |\mu|^4 dA.
\end{equation}

\subsubsection{Weil-Petersson sectional curvatures.}\label{subsec:sec-curv}
We now describe a lower bound on sectional curvatures which follows from \cite{Wol86}. We let $\mu_i, \mu_j \in T_X\sM_g^n$ be two orthogonal tangent vectors with $||\mu_i||_{WP}=||\mu_j||_{WP}=1$. We let $K^{WP}(\mu_i,\mu_j)$ be the \wep sectional curvature of the plane spanned by the real vectors corresponding to $\mu_i$ and $\mu_j$.
In \cite[Theorem 4.5]{Wol86}, Wolpert makes the following observations. Wolpert shows that
\beqar
\int_X D(|\mu_i||\mu_j|)|\mu_i||\mu_j|dA &\leq& \int_X D(|\mu_i|^2)|\mu_j|^2dA,\\
\left|\int_X D(\mu_i \mu_{\overline{j}})\mu_i \mu_{\overline{j}}dA\right|&\leq& \int_X D(|\mu_i||\mu_j|)|\mu_i||\mu_j|dA,\\
\left|\int_X D(\mu_i \mu_{\overline{j}}) \mu_{\overline{i}}\mu_jdA\right| &\leq& \int_X D(|\mu_i||\mu_j|)|\mu_i||\mu_j|dA.
\eeqar

Therefore as sectional curvature is given by 
\beqar
K^{WP}(\mu_i,\mu_j)&=&\Re{\int_X D(\mu_i \mu_{\overline{j}})\mu_i \mu_{\overline{j}}dA}-\frac{1}{2}\int_X D(\mu_i \mu_{\overline{j}})\mu_{\overline{i}}\mu_jdA\\
&&-\frac{1}{2}\int_X D(|\mu_i|^2)|\mu_j|^2dA
\eeqar
 putting these equations together gives
\begin{equation}\label{SecK-eq}
K^{WP}(\mu_i, \mu_j)\geq -2\int_X D(|\mu_i|^2)|\mu_j|^2dA.
\end{equation}

\subsubsection{Weil-Petersson Ricci curvatures.}\label{subsec:ric-curv}
Let $\{\mu_{i}\}_{i=1}^{3g-3+n}$ be a holomorphic orthonormal basis of $T_X\sM_g^n$. Then the Ricci curvature $\Ric(\mu_i)$ of $\sM^n_g$ at $X$ in the direction $\mu_{i}$ is given by
\begin{eqnarray*}
&&\Ric(\mu_{i})=-\sum_{j=1}^{3g-3+n}R_{i\overline{j}j\overline{i}}\\
&=&-\sum_{j=1}^{3g-3+n}\left(\int_{X} D(\mu_{i}\mu_{\overline{j}})\cdot (\mu_{j}\mu_{\overline{i}}) dA+\int_{X} D(|\mu_{i}|^2)\cdot (|\mu_{j}|^2) dA\right).
\end{eqnarray*}
Since $\int_X D(f)\cdot \overline{f} dA \geq 0$ for any function $f$ on $X$, by applying the argument in the proof of \eqref{SecK-eq} we have
\begin{equation}\label{RicK-eq}
-2\leq \frac{\Ric(\mu_i)}{\sum_{j=1}^{3g-3+n}\int_{X} D(|\mu_{i}|^2)\cdot (|\mu_{j}|^2)dA} \leq -1.
\end{equation}

\subsubsection{\wep Scalar Curvature.}
The scalar curvature $\Sca(X)$ at $X\in \sM_g^n$ is the trace of the Ricci tensor. We can express the scalar curvature as  
\begin{equation} \label{eqn:scalar curvature} 
\Sca(X)=-\sum_{i=1}^{3g-3+n}\sum_{j=1}^{3g-3+n}(\int_{X} D(\mu_{i}\mu_{\overline{j}})\cdot (\mu_{j}\mu_{\overline{i}}) dA+\int_{X} D(|\mu_{i}|^2)\cdot (|\mu_{j}|^2) dA).
\end{equation}
It is known from \cite[Proposition 2.5]{Wolf-W-1} that $-\Sca(X)$ is uniformly comparable to the quantity $|| \sum_{i=1}^{3g-3+n}|\mu_{i}|^2||_{WP}^2$. More precisely,
\begin{equation}\label{ScaK-eq}
 -2\int_{X}{(\sum_{i=1}^{3g-3+n}|\mu_{i}|^2)^2}dA \leq \Sca(X)\leq -\frac{1}{3}\int_{X}{(\sum_{i=1}^{3g-3+n}|\mu_{i}|^2)^2dA}.
\end{equation}


\section{Bounding the pointwise norm by the $L^2$ norm}\label{ub-hbd}
In this section we will bound the pointwise norm of a harmonic Beltrami differential $\mu = \overline\phi/\sigma$ in terms of its Weil-Petersson norm and the injectivity radius function. Our results will improve on prior work of Teo \cite{Teo09} and Wolpert \cite{Wolpert5}, giving the optimal asymptotics of Wolpert with its uniformity of Teo. As in Wolpert \cite[Proposition 7]{Wolpert5}, our approach will be to first decompose $\phi$ in the thin part of the surface into the leading and non-leading parts of its Laurent expansion. Then by a detailed analysis, we describe the leading term and give an explicit exponentially decaying upper bound on the non-leading term. 

Given $X\in \sM_g^n$ a hyperbolic surface of finite volume, for $z\in X$ we will let $r(z) = \inj(z)$ be the injectivity radius at $z$. We will refer several times to the a function $C(r)$ introduced by Teo in \cite{Teo09} which is given by
\bear \label{defn:C-inj}
C(r)&=&\left(\frac{4\pi}{3}\left(1-\sech^6\left(\frac{r}{2}\right)\right) \right)^{-\frac{1}{2}}\\
&=&\left(\frac{4\pi}{3}\left(1-\left(\frac{4e^{r}}{(1+e^{r})^2}\right)^3\right)\right)^{-\frac{1}{2}}. \nonumber \
\eear

It follows that $C(r)$ is decreasing with respect to $r$ and as $r$ tends to zero we have
$$C(r) = \frac{1}{\sqrt{\pi r}} + O(1) .$$ Furthermore $C(r)$ tends to $\sqrt{\frac{3}{4\pi}}$ as $r$ tends to infinity. 

Let $X=(S_g^n, \sigma(z)|dz|^2)\in \sM^n_g$ and $\phi \in Q(X)$ where $Q(X)$ is the space of  holomorphic quadratic differentials on $X$. We set
\bear
||\phi(z)||:=\frac{|\phi(z)|}{\sigma(z)} \quad \emph{for all $z\in X$},
\eear
and
\bear
||\phi||_2:= (\int_X ||\phi(z)||^2 \cdot \sigma(z)|dz|^2)^{\frac{1}{2}}.
\eear

We have the following result of Teo.

\begin{lemma}{(Teo, \cite[Proposition 3.1]{Teo09})} \label{lnubw}
Let $\phi \in Q(X)$ be a holomorphic quadratic differential on a  hyperbolic surface $X\in \sM^n_g$, and $r:X \rightarrow \R_+$ be the injectivity radius function. Then
 $$||\phi(z)||  \leq C(r(z))\cdot||\phi ||_2 =\frac{||\phi ||_2}{\sqrt{\pi} \cdot r(z)}(1+ o(r(z)))$$
where the constant $C(\cdot)$ is given by \eqref{defn:C-inj}.
\label{teo}
\end{lemma}

In \cite{Wolpert5}, Wolpert gave the following asymptotically optimal bound.

\begin{lemma}{(Wolpert, \cite[Corollary 11]{Wolpert5})}
Let $S$ be a surface of genus $g$ with $n$ punctures, and $X \in \mathcal M_g^n$ be any hyperbolic surface. Then for any $\epsilon > 0$ there exists a $\delta(\epsilon,S) > 0$  such that  if $\lsys \leq \delta(\epsilon,S)$ then for any $\phi \in Q(X)$ and $z\in X$   
$$\|\phi(z)\| \leq (1+\epsilon)\sqrt{\frac{2}{\pi}}\frac{\|\phi\|_2}{\sqrt{\lsys}}.$$
\label{wolpert_inj}
\end{lemma}

We will now derive a uniform bound that gives the asymptotics of Wolpert's bound above.

\subsection{Collar Neighborhoods}
We let $\phi \in Q(X)$ be a holomorphic quadratic differential on a Riemann surface $X\in \sM^n_g$ and $\gamma$ be  a simple closed geodesic of length $L$ in $X$. We lift  $\phi$ to $\tilde\phi$ on the  annulus  $A = \{ z\ | e^{-\frac{\pi^2}{L}} < |z| <  e^{\frac{\pi^2}{L}}\}$. Then $\tilde\phi(z) =\frac{f(z)}{z^2}dz^2$ where $f$ is holomorphic on $A$. Therefore we have the Laurent series
$$f(z) = \sum_{n=-\infty}^\infty a_nz^n.$$
We define
$$f_-(z) = \sum_{n<0} a_nz^n,\quad f_0(z) = a_0, \quad f_+(z) = \sum_{n>0} a_nz^n.$$
We therefore have the decomposition

$$\tilde \phi(z) = (f_-(z) +f_0(z)+f_+(z))\frac{dz^2}{z^2} =  \phi_-(z)+\phi_0(z)+\phi_+(z)$$

Let $\gamma \subset X\in \sM^n_g$ be a closed geodesic of length $L\leq2\eo$. By the Collar lemma (see \cite[Chapter 4]{Buser10}) there is an embedded collar $\mathcal{C}_\gamma$ of $\gamma$ in $X$ as follows.
\bear \label{def-collar}
\mathcal{C}_\gamma:= \{z\in X| \ d(z, \gamma) \leq \arcsinh(\frac{1}{\sinh(\frac{L}{2})})\}. 
\eear

We set 
\bear
||\phi|_{\Cc}||_2:= \left(\int_{\Cc} ||\phi(z)||^2 \cdot \sigma(z)|dz|^2\right)^{\frac{1}{2}}.
\eear

As $\mathcal C_\gamma$ embeds in $A$, we have that the injectivity radius function $r$ on $A$ coincides with the injectivity radius function on $\mathcal C_\gamma \subseteq X$. Also  if $z \in A$ has distance $d(z,\gamma)$ from the core closed geodesic then
$$\sinh(r(z)) =  \sinh(L/2)\cosh(d(z,\gamma))$$
Therefore it follows that 
$$\mathcal{C}_\gamma = \left\{z\in A \ | \ r(z) \leq \sinh^{-1}(\cosh(L/2))\right\}.$$
For $0<t \leq \sinh^{-1}(\cosh(L/2))$ we then define
$$C_t = \left\{z\in  A \ | \ r(z) \leq t\right\}.$$

In part of the following Proposition we will need to restrict to a sub-collar of the standard collar $\mathcal C_\gamma$. For this we define the constant $$\eoh = \frac{\log(3)}{2} = \sinh^{-1}\left(\frac{1}{\sqrt{3}}\right).$$ We prove the following

\begin{proposition}\label{fourier1}
Let $\phi \in Q(X)$ and $\mathcal C_\gamma$ be the collar about a closed geodesic $\gamma$ of length $L \leq 2 \eo$. Then 
\begin{enumerate}
 \item for any $z \in \mathcal C_\gamma$
 $$||\phi_0(z)|| \leq  \frac{1}{\sqrt{Lc_0(L)}}\frac{\sinh^2(L/2)}{\sinh^2(r(z))}\|\phi |_{\mathcal C_\gamma}\|_2 $$
 where
 \[c_0(L) = \cos^{-1}\left(\tanh\left(L/2\right)\right) + \frac{1}{2}\sin\left(2\cos^{-1}\left(\tanh\left(L/2\right)\right)\right) = \frac{\pi}{2} -\frac{L^3}{12}+O(L^5).\]

\item On $C_t$,  $\|\phi_\pm(z)\|$ attains its maximum on $\partial C_t$.

\item For $z \in \mathcal C_\gamma$ in the sub-collar $C_{\overline{\epsilon}_2}=\{z\in A\ |\ r(z) \leq  \eoh\}$,
$$||\phi_{\pm}(z)|| \leq F(r(z))\|\phi |_{\mathcal C_\gamma}\|_2$$
 where
$$F(r(z)) =\frac{e^{\pi\sqrt{3}} C\left(\eoh\right) e^{-\frac{\pi}{\sinh(r(z))}}}{3\sinh^2(r(z))}\leq C\left(\eoh\right).$$
 \item For $z \in \mathcal C_\gamma$ in the sub-collar $C_{\overline{\epsilon}_2}$ with $r(z) \leq \eoh$ 
 $$||\phi(z)|| \leq G(r(z))\|\phi |_{\mathcal C_\gamma}\|_2 $$
 where 
 $$G(r) =  \frac{1}{\sqrt{2rc_0(2r)}}+\frac{2e^{\pi\sqrt{3}} C(\eoh)e^{-\frac{\pi}{\sinh(r)}}}{3\sinh^2(r)} = \frac{1}{\sqrt{\pi r}}\left(1+\frac{2r^3}{3\pi} + O(r^5)\right).$$
\item For $z \in \mathcal C_\gamma$ with $r(z) \leq \eo$ then
 $$||\phi(z)|| \leq \frac{||\phi ||_2}{\sqrt{r(z)}}.$$
\end{enumerate}
\label{decomp}
\end{proposition}
\bp
Let  $S = \{ z = x+iy \ | \ |y| < \pi/2\}$ be the strip, then the the hyperbolic metric on $S$ is $\rho_S(z) = |dz|/\cos(y)$. By the Collar Lemma \cite[Theorem 4.1.6]{Buser10} the injectivity radius function  on $S$ satisfies
\bear\label{dist-form}
\sinh(r(z)) = \frac{\sinh(L/2)}{\cos(y)}.
\eear

\noindent We have the $\Z$ cover $\pi:S \rightarrow A$ given by $\pi(z) = e^{\frac{2\pi iz}{L}}$.
Therefore the hyperbolic metric on $A$ is given by 
$$\rho(z) = \frac{L}{2\pi} \frac{1}{|z|\cos\left(\frac{L}{2\pi} \log|z|\right)}.$$
It follows that $\mathcal C_\gamma$ lifts to the strip $\mathcal S_\gamma =  \{ w = x+iy \ | \ |y| < h(L)\}$ where $$h(L) = \cos^{-1}(\tanh(L/2)).$$  Therefore $\mathcal C_\gamma =\{z\in A\ | \ e^{-s(L)} < |z| < e^{s(L)}\}$ where $$s(L) = 2\pi\cdot \frac{h(L)}{L}.$$ 

We first show that $\phi_-,\phi_0,\phi_+$ are all orthogonal on $\mathcal C_\gamma$. We have 
\beqar
\|\phi |_{\mathcal C_\gamma}\|_2^2 &=& \int_{\mathcal C_\gamma} \frac{|\phi(z)|^2}{\rho^2(z)} = \sum_{n,m}\int_{e^{-s(L)}}^{e^{s(L)}}\int_0^{2\pi} \frac{a_n\overline a_m z^n \overline{z}^m}{|z|^4 \rho^2(r)} rdrd\theta \\
&=& \sum_{n,m}\left(\int_{e^{-s(L)}}^{e^{s(L)}} \frac{a_n\overline a_m r^{n+m-3}}{ \rho^2(r)} dr\right)\left(\int_0^{2\pi} e^{i(n-m)\theta}d\theta\right) \\
&=& 2\pi\sum_{n}\int_{e^{-s(L)}}^{e^{s(L)}} \frac{|a_n|^2 r^{2n-3}}{ \rho^2(r)} dr.
\eeqar

\noindent Therefore 
\begin{equation}
\|\phi |_{\mathcal C_\gamma}\|_2^2 =  \|\phi_-|_{\mathcal C_\gamma}\|_2^2+\|\phi_0 |_{\mathcal C_\gamma}\|_2^2+\|\phi_+ |_{\mathcal C_\gamma}\|_2^2.
\label{parseval}
\end{equation}
This gives the bound
$$\|\phi |_{\mathcal C_\gamma}\|_2^2  \geq  \|\phi_0 |_{\mathcal C_\gamma}\|_2^2 = 2\pi |a_0|^2 \frac{4\pi^2}{L^2}\int_{e^{-s(L)}}^{e^{s(L)}} \frac{\cos^2\left(\frac{L}{2\pi} \log r\right)}{r} dr.$$
We let $t = \frac{L}{2\pi} \log r$ giving $dt = \frac{L}{2\pi r} dr$ and
$$\|\phi_0 |_{\mathcal C_\gamma}\|_2^2 =  |a_0|^2 \frac{16\pi^4}{L^3}\int_{-h(L)}^{h(L)} \cos^2(t)dt.
 $$
We define
 $$c_0(L) = \int_{-h(L)}^{h(L)} \cos^2(t) dt  =  h(L)+ \frac{1}{2}\sin(2h(L)) = \frac{\pi}{2} -\frac{L^3}{12}+ O(L^5).$$
 Then
 $$\|\phi_0 |_{\mathcal C_\gamma}\|_2^2 =  |a_0|^2 \frac{16\pi^4}{L^3}c_0(L).$$
For $z \in \mathcal C_\gamma$, we have 
 $$\|\phi_0(z)\| = \frac{4\pi^2|a_0|}{L^2} \cos^2\left(\frac{L}{2\pi} \log|z|\right) =\frac{1}{\sqrt{Lc_0(L)}}\frac{\sinh^2(L/2)}{\sinh^2(r(z))}\|\phi_0 |_{\mathcal C_\gamma}\|_2$$
where in the last equality we apply the following version of formula \eqref{dist-form}
\[\cos\left(\frac{L}{2\pi} \log|z|\right) = \frac{\sinh(L/2)}{\sinh(r(z))}.\]
Thus
\bear
\|\phi_0(z)\|  \leq \frac{1}{\sqrt{Lc_0(L)}}\frac{\sinh^2(L/2)}{\sinh^2(r(z))}\|\phi |_{\mathcal C_\gamma}\|_2
\eear
giving (1).

We consider $\phi_+(z) = f_+(z)dz^2/z^2$. We have that $f_+(z)$ is holomorphic on the disk $D_+  = \{ z\ |  \ |z| <  e^{\frac{\pi^2}{L}}\}$. Furthermore $f_+(z)/z$ extends holomorphically to $D_+$. By the maximum principle  the maximum modulus of $f_+(z)/z$ on $B(s) = \{ z\ | \ |z| \leq s\}$ is on the boundary. Therefore the maximum modulus of $f_+(z)/z$ on $B(s)$ is at some $z_s \in \partial B(s)$ with $M_s = |f(z_s)|/|z_s|$. We have for $z \in B(s)$ 
$$||\phi_+(z)|| = \frac{|f_+(z)|}{|z|^2}\cdot\frac{4\pi^2}{L^2}|z|^2\cos^2\left(\frac{L}{2\pi} \log|z|\right) \leq M_s\frac{4\pi^2}{L^2}|z|\cos^2\left(\frac{L}{2\pi} \log|z|\right).$$
Recall that
\[||\phi_+(z_s)||= M_s \frac{4\pi^2}{L^2} s\cos^2(\left(\frac{L}{2\pi} \log s\right).\]
Therefore
\bear\label{upp-+}
||\phi_+(z)|| \leq \frac{||\phi_+(z_s)||}{s\cos^2(\left(\frac{L}{2\pi} \log s\right)}\cdot \left(|z|\cos^2\left(\frac{L}{2\pi} \log|z|\right)\right).
\eear
We observe that $x\cos^2\left(\frac{L}{2\pi} \log x\right)$ is monotonically increasing on $[1, e^{s(L)}]$. To see this, we consider equivalently the function $u(t) = e^{2\pi t/L} \cos^2(t)$ on $[-h(L),h(L)]$. Differentiating it we get
$$u'(t) = 2e^{2\pi t/L}\cos(t)\left(\frac{\pi}{L}\cos(t)-\sin(t)\right).$$
Thus $u$ is monotonic for  $\tan(t) \leq \frac{\pi}{L}.$ As $t \leq h(L) = \cos^{-1}(\tanh(L/2))$ we have 
$$\tan(t) \leq \tan(h(L)) = \frac{1}{\sinh(L/2)} \leq \frac{2}{L} \leq \frac{\pi}{L}.$$  
Thus $u$ is monotonic on $[1, \frac{L}{2\pi}\cdot s(L)]$.
Therefore  $\|\phi_+(z)\|$ has maximum modulus in $C_t$ on the boundary. Similarly one may prove that $\|\phi_-(z)\|$ has maximum modulus in $C_t$ on the boundary by using $\frac{1}{z}$ as a variable. This proves (2).

To prove (3) we use Teo's bound from Lemma \ref{teo}. By Teo
$$\|\phi_+(z_s)\| \leq C(r(z_s))\cdot||\phi_+ |_{B(z_s,r(z_s))} ||_2$$
where $B(z,r)$ is the hyperbolic ball about $z$ of radius $r$. We choose $z_s$ in the collar such that $B(z_s,r(z_s)) \subseteq \mathcal C_\gamma$. By the Collar Lemma \cite[Theorem 4.1.6]{Buser10}, a point of injectivity radius $r$ is a distance $d$ from the boundary of the collar where
$$\sinh(r) = \cosh(\frac{L}{2})\cosh d - \sinh d.$$
We note that solving $d = r$ gives 
$$r = \tanh^{-1}\left(\frac{\cosh(\frac{L}{2})}{2}\right) \geq \tanh^{-1}(1/2).$$
Therefore we choose $z_s$ such that $r(z_s) =  \tanh^{-1}(1/2) = \eoh.$ Then by Lemma \ref{teo} and \eqref{parseval}
$$\|\phi_+(z_s)\| \leq C(\eoh)\cdot||\phi_+ |_{\mathcal C_\gamma} ||_2 \leq C(\eoh)\cdot||\phi|_{\mathcal C_\gamma} ||_2 .$$
This together with \eqref{upp-+} implies that
$$||\phi_+(z)|| \leq \frac{C(\eoh)}{s\cos^2(\left(\frac{L}{2\pi} \log s\right)}\cdot ||\phi |_{\mathcal C_\gamma}  ||_2 \left(|z|\cos^2\left(\frac{L}{2\pi} \log|z|\right)\right).$$
Recall that \eqref{dist-form} gives
$$\cos\left(\frac{L}{2\pi} \log|z|\right) = \frac{\sinh(L/2)}{\sinh(r(z))}.$$
Therefore
$$|z| = e^{\pm\frac{2\pi}{L}\left(\cos^{-1}\left(\frac{\sinh(L/2)}{\sinh(r(z))}\right)\right)}$$
where the sign depends on which side of the core closed geodesic you are on. 
We  rewrite the bound in terms of injectivity radius. Recall that $s>1$. Then for 
$|z| \geq 1$, i.e., $|z| = e^{\frac{2\pi}{L}\left(\cos^{-1}\left(\frac{\sinh(L/2)}{\sinh(r(z))}\right)\right)}$, 
$$||\phi_+(z)|| \leq \frac{C(\eoh)\sinh^2(\eoh)e^{\frac{2\pi}{L}\left(\cos^{-1}\left(\frac{\sinh(L/2)}{\sinh(r(z))}\right)-\cos^{-1}\left(\frac{\sinh(L/2)}{\sinh(\eoh)}\right)\right)}}{\sinh^2(r(z))}\cdot||\phi |_{\mathcal C_\gamma}  ||_2.$$
Note that $\sinh(\eoh) = 1/\sqrt{3}$. Also  for $0< x < y <\pi  $ then $x-y\leq \cos(y)-\cos(x)$ giving
$$||\phi_+(z)||^2 \leq \frac{C(\eoh)e^{-\frac{2\pi\sinh(L/2)}{L}\left(\frac{1}{\sinh(r(z))}-\sqrt{3}\right)}}{3\sinh^2(r(z))}\cdot||\phi |_{\mathcal C_\gamma}  ||_2.$$
As $\sinh(x) \geq x$ we have for $|z|\geq 1$,
\bear \label{3.10}
\ \ \ \ \ \ \ \  ||\phi_+(z)|| \leq \frac{C(\eoh)e^{-\pi\left(\frac{1}{\sinh(r(z))}- \sqrt{3}\right)}}{3\sinh^2(r(z))}\cdot||\phi |_{\mathcal C_\gamma}  ||_2 = F(r(z))\cdot||\phi |_{\mathcal C_\gamma}  ||_2.
\eear

\noindent We note that $r(z) = r(1/z)$. Also by the above, the maximum of $\|\phi_+(z)\|$ on $\{z\ |\ 1/c \leq |z| \leq c\}$ is on the boundary $|z| = c$ where $1<c\leq e^{s(L)}$. Therefore for $|z|= 1/c < 1$ we have $$||\phi_+(z)|| \leq\max_{|w| = 1/c}||\phi_+(w)|| \leq \max_{|w| = c}||\phi_+(w)|| \leq F(r(z))\cdot||\phi |_{\mathcal C_\gamma}  ||_2.$$

\noindent Thus for $r(z) \leq \eoh$
\bear
\quad ||\phi_+(z)|| \leq \left(\frac{C(\eoh)e^{\pi\sqrt{3}} e^{-\frac{\pi}{\sinh(r(z))}}}{3\sinh^2(r(z))}\right)||\phi |_{\mathcal C_\gamma}  ||_2 \leq C(\eoh)||\phi |_{\mathcal C_\gamma}  ||_2
\eear
where in the last inequality we apply that $\frac{e^{-\frac{\pi}{\sinh(r(z))}}}{\sinh^2(r(z))}$ is increasing. Similar as in the proof of Part (2) if we consider $\frac{1}{z}$ as a variable, one may also get the same bound for $||\phi_-(z)||$. This proves (3).

For proving (4), we combine the bounds above using
$$\|\phi(z)\| \leq \|\phi_-(z)\|+\|\phi_0(z)\|+\|\phi_+(z)\|.$$
First observe that both $\frac{\sinh{(L/2)}}{\sqrt{L}}$ and $\frac{\sinh{(L/2)}}{\sqrt{c_0(L)}}$ are increasing. Since $2r(z)\geq L$,
for any $z \in \mathcal C_\gamma$ we have
 $$||\phi_0(z)|| \leq  \frac{1}{\sqrt{2r(z)c_0(2r(z))}}||\phi |_{\mathcal C_\gamma}\|_2.$$

\noindent Therefore for $z\in \mathcal C_\gamma$ with $r(z)\leq \overline{\epsilon}_2$,
\bear
\|\phi(z)\| \leq G(r(z))||\phi |_{\mathcal C_\gamma}  ||_2
\eear
where
$$G(r) = \frac{1}{ \sqrt{2rc_0(2r)}}+\frac{2e^{\pi\sqrt{3}} C(\eoh)e^{-\frac{\pi}{\sinh(r)}}}{3\sinh^2(r)}.$$
This proves (4).

To prove (5) we combine the above bound  for $r(z) \leq \eoh$ with Teo's bound for $r(z) \leq \eo$. If $r(z) \geq \eoh$ by Lemma \ref{teo} we have that
$$\|\phi(z)\| \leq   C(r(z))\cdot ||\phi  ||_2 \leq \sqrt{r(z)}C(r(z))\cdot \frac{||\phi  ||_2}{\sqrt{r(z)}}.$$
As $C(x)\sqrt{x}$ is monotonically decreasing with $C(\eoh)\cdot \sqrt{\eoh} = .8091$ we have
$$\|\phi(z)\| \leq    \sqrt{\eoh}C(\eoh) \cdot \frac{||\phi  ||_2}{\sqrt{r(z)}} =.8091 \frac{||\phi  ||_2}{\sqrt{r(z)}}.$$

\begin{figure}[htbp] 
   \centering
   \includegraphics[width=3in]{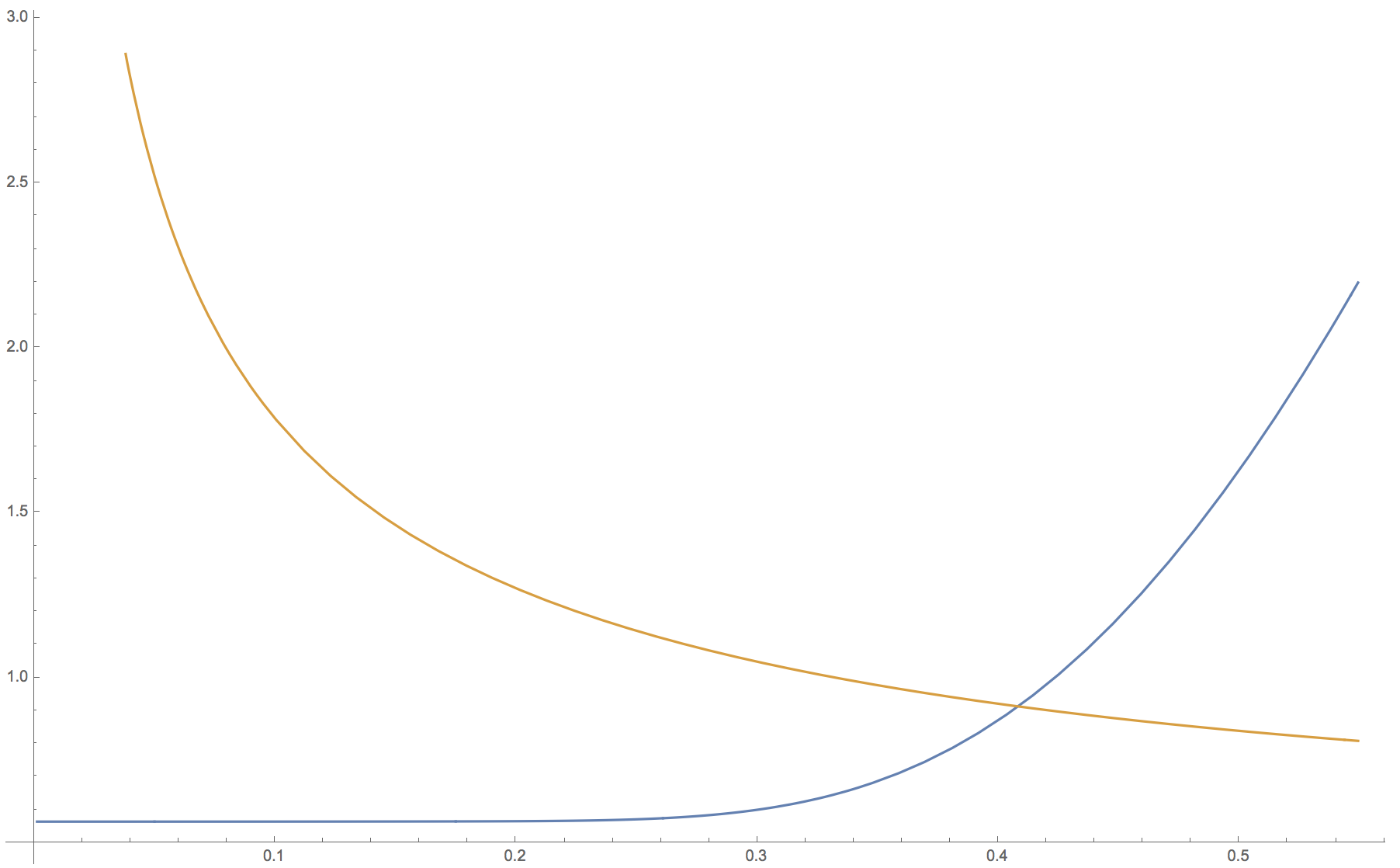} 
   \caption{Plot of $H(r)$ and $\sqrt{r}C(r)$ on $(0,\eoh]$}
   \label{m}
\end{figure}

\noindent We now consider $r(z) \leq \eoh$. We have that $H(r) = G(r) \cdot \sqrt{r}$ is monotonically decreasing. Therefore Part (4) above together with Lemma \ref{teo} imply that 
$$\|\phi(z)\| \leq   \min\left(H(r(z)), \sqrt{r(z)}C(r(z)\right) \cdot \frac{||\phi  ||_2}{\sqrt{r(z)}}.$$
Considering $m(r) = \min(H(r), \sqrt{r}C(r))$ on $(0,\eoh]$ we have by computation that $m(r) \leq  m_0 = .9137$ (see figure \ref{m}). 

Therefore  for $r(z) \leq \eo$
\begin{equation}
 \|\phi(z)\| \leq  \max\{.8091,\  .9137\}\cdot \frac{||\phi  ||_2}{\sqrt{r(z)}} \leq \frac{||\phi  ||_2}{\sqrt{r(z)}}
\label{min_r_bound}
\end{equation}
which completes the proof.
\ep

\subsection{Cusp neighborhoods}
We now consider the cusp neighborhoods of $X \in \sM^n_g$. Then each cusp $c$ gives a cover $\pi:\Delta^*\rightarrow X$ where $\Delta^*=\{ z\ | \ 0< |z| < 1\}$. The hyperbolic metric on $\Delta^*$ is $\rho(z) = -1/|z|\log|z|$.

By the Collar Lemma (see \cite[Chapter 4]{Buser10}), $c$ has a collar $\mathcal C_c$ which lifts to $A_c = \{ z\ | \ 0 < |z| < e^{-\pi}\}$ with $\pi$ injective on $A_c$. Furthermore  as $\mathcal C_c$ is embedded, the  injectivity radius function $r$ on $X$ lifts to the injectivity radius function on $A_c$ with  $A_c = \{ z\in A \ | \ r(z) < \eo\}.$
We have
\begin{lemma}\label{up-cusp}
Let $X \in \sM^n_g$ and $\phi \in Q(X)$. If $z \in \mathcal C_c$, then
$$\|\phi(z)\| \leq K(r(z))\|\phi\|_2\leq C(\eo)||\phi||_2$$
where 
$$K(r) = \left(\frac{C(\eo)e^{\pi}e^{-\frac{\pi}{\sinh(r)}}}{\sinh^2(r)}\right)$$
and $C(\eo)=.7439$.
\end{lemma}
\bp
As before we have $\phi = \phi_- + \phi_0 +\phi_+$. We have the hyperbolic metric on $A_c$ is $\rho(z) = -1/|z|\log|z|$. The lemma is trivially true if $||\phi||_2 = \infty$. Therefore we consider $||\phi||_2<\infty$.  It follows that $\phi_0 = \phi_- = 0$. We now bound $\|\phi(z)\|$ as above. If $\phi(z) = f(z)dz^2/z^2$ then $f(z)/z$ extends to $B(s) = \{ z\ | \  |z| < s\}$ and has maximum modulus at $z_s$ with $|z_s| = s$. Therefore
$$||\phi(z)|| \leq \frac{||\phi(z_s)||}{s(\log|s|)^2}\cdot|z| (\log|z|)^2.$$
It can easily be checked that $|z| (\log|z|)^2$ is monotonic on $A_c$. By the Collar Lemma, the  the injectivity radius on $A_c$ satisfies
$\sinh(r(z)) = -\pi/ \log|z|$. Therefore by letting $s = e^{-\pi}$ (the maximal cusp) and using Lemma \ref{teo} we obtain that for $r(z) \leq \eo$,
$$||\phi(z)|| \leq  \left(\frac{C(\eo)e^{\pi}e^{-\frac{\pi}{\sinh(r(z))}}}{\sinh^2(r(z))}\right)||\phi||_2 = K(r(z))\cdot||\phi||_2.$$
The function $\frac{e^{-\frac{\pi}{\sinh(r)}}}{\sinh^2(r)}$ is monotonically increasing on $[0, \epsilon_2]$. Recall that $\sinh{(\eo)}=1$. So we have
\[||\phi(z)||\leq C(\eo)||\phi||_2.\]
Which completes the proof.
\ep

\subsection{Uniform upper bounds for $||\phi||$} In this subsection we discuss several applications of Proposition \ref{fourier1} and Lemma \ref{up-cusp}. The first one is to show Proposition \ref{lnubw-0}.
\bp [Proof of Proposition \ref{lnubw-0}]
Let $z\in X$ with $\inj(z)\leq \eo$. Then $z$ is in either a collar or a cusp. If $z$ is in a collar, the claim follows by Part (5) of Proposition \ref{fourier1}. If $z$ is in a cusp, the claim follows by Lemma \ref{up-cusp}.  
\ep

We define $\ell^+_{sys}(X) = \min(2\eo, \lsys)$. Then we have

\begin{corollary}
Let $X \in \sM^n_g$  and $\phi \in Q(X)$. Then $$ \|\phi\|_\infty    \leq \sqrt{\frac{2}{\ell^+_{sys}(X)}} ||\phi||_2.$$
\end{corollary}
\label{systole_bound}
\bp
If $r(z) \geq \eo$  or $z$ is in a cusp neighborhood then as $\ell^+_{sys}(X) \leq2\eo$, it follows by Lemma \ref{teo} and Lemma \ref{up-cusp} that
$$||\phi(z)||\leq C(\eo)||\phi||_2 \leq \sqrt{2\eo}C(\eo)\frac{||\phi||_2}{\sqrt{\ell^+_{sys}(X)}}.$$
We have $\sqrt{2\eo}\cdot C(\eo) = .9877 < \sqrt{2}.$ So the claim follows for these two cases.

If $z$ is in a collar neighborhood with $r(z)\leq \eo$,  it follows by \eqref{min_r_bound} that
\begin{equation}
||\phi(z)||\leq   m_0 \frac{||\phi||_2}{\sqrt{r(z)}} \leq \sqrt{2}\cdot m_0 \cdot \frac{||\phi||_2}{\sqrt{\lsys}}.
\label{min_sys_bound} 
\end{equation}
The claim also follows as $m_0  < .9137$.
\ep

\br
We note that we can use  Proposition \ref{decomp} to  give a bound for Wolpert's  Lemma \ref{wolpert_inj} which is independent of topology.
We let $H(r) = G(r)\cdot \sqrt{r}$. Then $H(r)$ is monotonically increasing with 
$$\lim_{r\rightarrow 0} H(r) = \frac{1}{\sqrt{\pi}}.$$
We note for from Part (4) of Proposition \ref{decomp} that for $r(z) \leq \eoh$  
\beqar
\|\phi(z)\| \leq  \frac{\|\phi\|_2}{\sqrt{r(z)}}.
\eeqar
Thus for $ \pi/2\cdot \lsys \leq r(z) \leq \eoh$ we have
\bear
\|\phi(z)\| \leq  \frac{\|\phi\|_2}{\sqrt{r(z)}} \leq \sqrt{\frac{2}{\pi}}\frac{\|\phi\|_2}{\sqrt{\lsys}}.
\eear
We choose $\delta_1$ such
$$\delta_1 = \frac{2}{\pi}H^{-1}\left(\frac{1+\epsilon}{\sqrt{\pi}}\right).$$
Then it follows by Part (4) of Proposition \ref{decomp} that for $\lsys < \delta_1$ and $r(z) \leq \min\{\frac{\pi}{2}\cdot \lsys, \eoh\} \leq \min\{ H^{-1}\left(\frac{1+\epsilon}{\sqrt{\pi}}\right), \eoh\}$
\bear
\|\phi(z)\| \leq (1+\epsilon)\sqrt{\frac{2}{\pi}}\frac{\|\phi\|_2}{\sqrt{\lsys}}.
\eear
Now for $r(z) \geq \eoh$  as $C(\eoh) = 1.09 < 2$
$$\|\phi(z)\| \leq C(\eoh)\|\phi\|_2 \leq 2\|\phi\|_2.$$
Thus for  $\lsys < \frac{1}{2\pi}$  and $r(z) \geq \eoh$ we have
\bear
\|\phi(z)\| \leq \sqrt{\frac{2}{\pi}}\frac{\|\phi\|_2}{\sqrt{\lsys}}.
\eear
We therefore choose $\delta = \min(\delta_1,\frac{1}{2\pi})$ to get the following result.
\bt
Let $X \in \mathcal \sM^n_g$ be any hyperbolic surface. Then for any $\epsilon > 0$ there exists a constant $\delta(\epsilon) > 0$ only depending on $\eps$ such that  if $\lsys \leq \delta(\epsilon)$ then for any $\phi \in Q(X)$ and $z\in X$,   
$$\|\phi(z)\| \leq (1+\epsilon)\sqrt{\frac{2}{\pi}}\frac{\|\phi\|_2}{\sqrt{\lsys}}.$$
\label{wolpert_inj}
\et

\noindent We note by the expansion of $G$ we have for $\epsilon$ small,
$$\delta(\epsilon) = \frac{2}{\pi}H^{-1}\left(\frac{1+\epsilon}{\sqrt{\pi}}\right) \simeq \left(\frac{12\epsilon}{\pi^2}\right)^{1/3}.$$
\er

\subsection{Fixing the length of short curves}
Let $X \in \mathcal M^n_g$ and for $\alpha$ a closed curve, we let $l_\alpha$ be the geodesic length function on $\mathcal M^n_g$. Then we let $dL_\alpha \in T^*(M^n_g)$ be the complex one-form such that $\Re dL_\alpha  = dl_\alpha$.   
We define 
\bear \label{def-P}
\emph{$P(X) \subseteq T_X^*(\mathcal M_g)$  = $span\{(dL_\alpha)_X\ | \ l_\alpha(X) \leq \eo\}$} 
\eear
and
\bear\label{def-P-p}
P(X)^{\perp} =  \{ \mu\ |\ \langle\phi,\mu\rangle = 0, \forall\ \phi \in P(X)\}  \subseteq T_X(\mathcal M^n_g).
\eear 

\noindent The plane $P(X)^{\perp}$ is the set of directions that fix the length of short curves. We have the following immediate consequence of Proposition \ref{fourier1}.

\begin{lemma}\label{ub-v-orth}
Let $\mu \in P(X)^{\perp}$ then
$$\|\mu(z)\|  \leq \sqrt{2}\cdot\|\mu\|_2.$$
Furthermore for $r(z) \leq \eoh$
$$\|\mu(z)\| \leq 2\cdot F(r(z))\cdot\|\mu\|_2.$$
Where $F(r(z))$ is defined in Proposition \ref{fourier1}.
\end{lemma}
\bp
Let $\mu = \overline{\phi}/\rho^2  \in P(X)^{\perp}$.
Recall that $C(\eoh)= 1.0917$. If $r(z) \geq \eoh$, then by Lemma \ref{teo}
\bear
\|\mu(z)\| \leq C(\eoh)\cdot \|\mu\|_2\leq \sqrt{2}\|\mu\|_2.
\eear
Similarly if $z$ is in a cusp neighborhood, then 
\bear
\|\mu(z)\| \leq K(r(z)) \leq C(\eo)\cdot \|\mu\|_2\leq C(\eoh)\|\mu\|_2.
\eear

\noindent Now we consider the remaining case. That is, $r(z) \leq \eoh$ and $z \in \mathcal C_\alpha$ where $\alpha\subset X$ is a closed geodesic with $l_\alpha(X) \leq 2\eoh$.
We lift $\phi$ to $\hat\phi$ on the annulus $A$ and have as before $\hat\phi(z) = \phi_-+\phi_0+\phi_+$ with  $\phi_0(z) = a\frac{dz^2}{z^2}$ for $a \in \C$. By the Gardiner formula \cite{IT92-book} we have
\beqar
0 &=& \langle dL_\alpha,\mu\rangle = \frac{2}{\pi}\int_A \frac{\overline{\hat\phi(z)}}{\rho(z)^2} \frac{dz^2}{z^2} \\
&=& \frac{2}{\pi}\int_A \frac{\overline{\phi_0(z)}}{\rho(z)^2} \frac{dz^2}{z^2} = \frac{2a}{\pi} \int_A \frac{dxdy}{r^4\rho^2(r)} = a l_\alpha(X).
\eeqar
  
\noindent Therefore $a = 0$ and $\phi_0 = 0$ (see also \cite[Proposition 8.5]{Wolf12}). Then it follows from Part (3) of Proposition \ref{fourier1} that
\bear \label{fl-eq-0-0}
\|\phi(z)\| \leq \|\phi_-(z)\| +\|\phi_+(z)\|  \leq 2\cdot F(r(z))\|\phi\|_2
\eear
where
\[F(r)=\frac{e^{\pi\sqrt{3}} C(\eoh)e^{-\frac{\pi}{\sinh(r)}}}{3\sinh^2(r)}.\]

\begin{figure}[htbp] 
   \centering
   \includegraphics[width=3in]{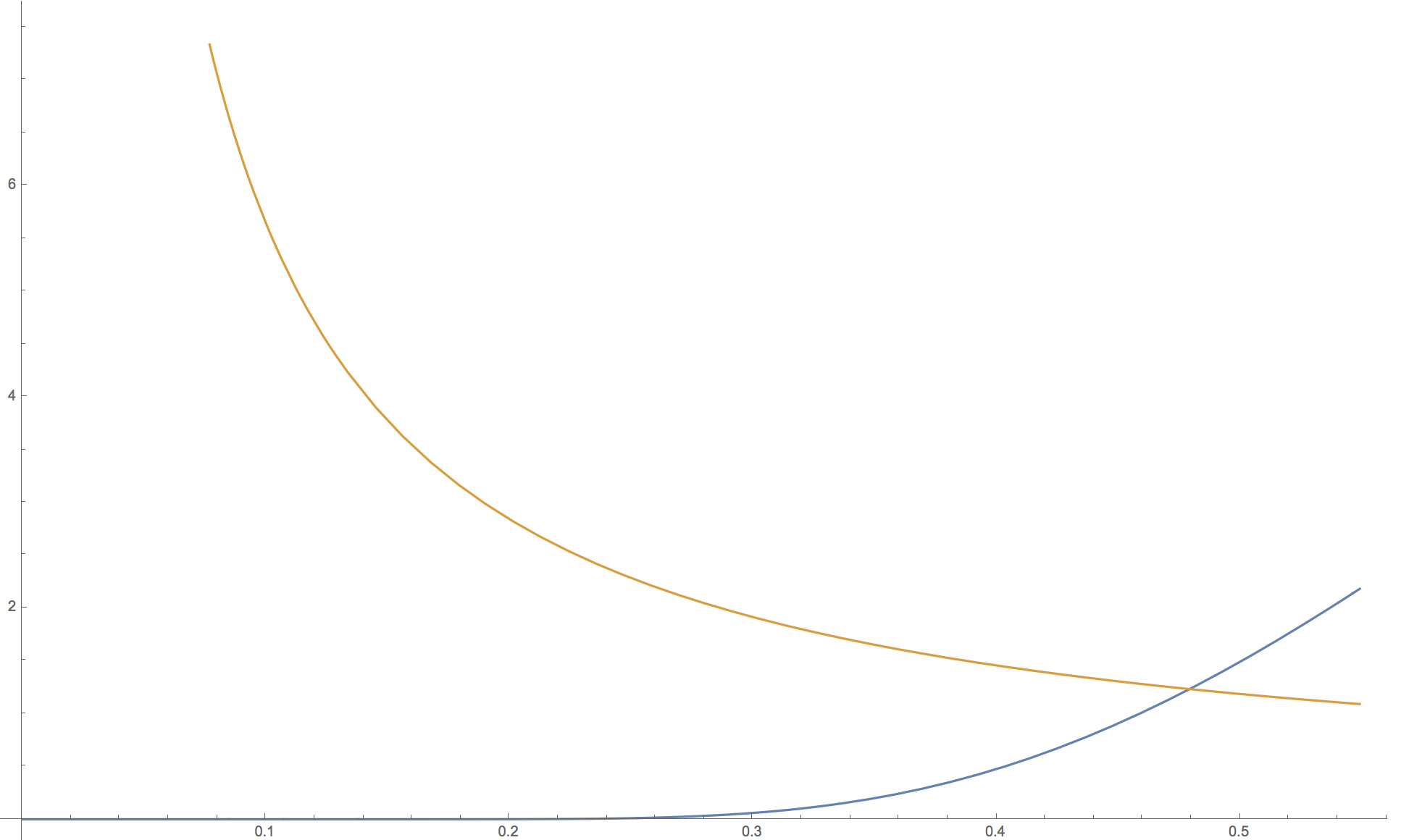} 
   \caption{Plot of $2F(r)$ and $C(r)$ on $(0,\eoh]$}
   \label{m'}
\end{figure}

\noindent Together with Lemma \ref{teo}, by letting $m'(r) =  \min(2F(r),C(r))$ we have
 $$\|\phi(z)\| \leq m'(r(z))\|\phi\|_2.$$
\noindent On $(0,\eoh]$ by computation we have $m'(r) \leq 1.2333$ 
(see figure \ref{m'}). 

Therefore
\bear
\|\phi(z)\| \leq  \sqrt{2}\cdot\|\phi\|_2
\eear
and proving the first inequality. 

We note that $K(r) \leq 2F(r)$ on $(0,\eoh]$ where $K(r)$ is defined in Lemma \ref{up-cusp}. Then it follows by Lemma \ref{up-cusp} and \eqref{fl-eq-0-0} for all $r(z) \leq \eoh$
$$\|\mu(z)\| \leq 2\cdot F(r(z))\cdot\|\mu\|_2$$
which completes the proof.
\ep

\section{Uniform lower bounds for \wep curvatures}\label{s-WP-curv}
The following bounds is essentially due to Teo \cite{Teo09}. As we need a slightly modified version, we give  the following version due to Ken Bromberg.
\begin{proposition}\label{Teo bound}
Fix $z\in X$ and let $U \subset T_X\sM_g^n$ be a subspace and $K_z>0$  a constant such that for all harmonic Beltrami differentials $\mu\in U$ we have
$$\|\mu(z)\| \le K_z \|\mu\|_2.$$
Then if $\mu_1, \dots, \mu_k$ is an orthonormal family in $U$ we have
$$\sum_{i=1}^{k}\|\mu_i(z)\|^2 \le K_z^2.$$
\end{proposition}

\bp
Pick constants $c_1, \dots, c_k$ such that $|c_i| = \|\mu_i(z)\|$ and the directions of maximal and minimal stretch of the Beltrami differentials $c_i\mu_i$ all agree at $z$.\footnote{For example if we choose a chart near $z$, in the chart the $\mu_i$ are realized by functions and we can let $c_i = \overline{\mu_i(z)}$. Then, in this chart, the directions of maximal and minimal stretch at $z$ of each $c_i \mu_i$ are the real and imaginary axis.} We then let $$\mu_z = \sum_{i=1}^k c_i\mu_i$$ and observe that our conditions on the directions of maximal and minimal stretch give that
$$\|\mu_z(z)\| = \sum_{i=1}^k \|c_i\mu_i(z)\| = \sum_{i=1}^{k} \|\mu_i(z)\|^2.$$
As the $\mu_i$ are orthonormal we also have
$$\|\mu_z\|^2 = \sum_{i=1}^{k} |c_i|^2=\sum_{i=1}^k \|\mu_i(z)\|^2.$$
As $\mu_z$ is a linear combination of harmonic Beltrami differentials it is also a harmonic Beltrami differential so 
$$\|\mu_z(z)\| \le K_z \|\mu_z\|$$
and therefore
$$\|\mu_z(z)\|^2 \le K_z^2\|\mu_z\|^2 = K_z^2 \|\mu_z(z)\|.$$
Dividing by $\|\mu_z(z)\| = \sum_{i=1}^k \|\mu_i(z)\|^2$ gives the result.
\ep

In this section we prove Theorem \ref{i-LUB-Ric}. Before proving it, we provide a uniform upper bound for any holomorphic orthonormal frame at $X \in \sM_g^n$.

First we make a thick-thin decomposition of $X\in\sM_g^n$ into three pieces as follows. Let $\eo$ be the Margulis constant as in previous sections. We set 
\[X_1:=\{q\in X:\ \inj(q)\geq \eo\},\]
\[X_2:=\{q\in \emph{cusps}:\ \inj(q)< \eo\},\]
\[X_3:=\{q\in \emph{collars}:\ \inj(q)< \eo\}.\]
So $X=\bigcup_{i=1}^3 X_i$. We note that the set $X_2$ and $X_3$ may be empty. Actually Buser and Sarnak \cite{BS94} showed that $\sup_{X \in \sM_g} \inj(X)\asymp \ln (g)$ for all $g\geq 2$.

Let $\{\mu_{i}\}_{i=1}^{3g-3+n}$ be a holomorphic orthonormal basis of $T_X\sM_g^n$. Our aim is to bound $\sum_{i=1}^{3g-3+n}|\mu_{i}|^2(z)$ from above.

First we restrict the discussion on $X_1$. In this case, Teo's formula \cite[Equation (3.12)]{Teo09}, which extends to the punctured case by Proposition \ref{Teo bound}, gives
\be \label{ub-l-1/5}
\sup_{z \in X_1}\sum_{i=1}^{3g-3+n}|\mu_{i}|^2(z)\leq C(\eo)^2= .5533.
\ene
This bound is an easy application of Lemma \ref{lnubw} and Proposition \ref{Teo bound}.

Next we consider the case on $X_2$. Recall that Lemma \ref{up-cusp} says that for any $x \in X_2$, $||\phi(z)||\leq C(\eo)||\phi||_2$. Therefore it follows by Proposition \ref{Teo bound} that
\be \label{ub-l-1/5-cusp}
\sup_{z \in X_2}\sum_{i=1}^{3g-3+n}|\mu_{i}|^2(z)\leq C(\eo)^2= .5533
\ene

Now we deal with the case on $X_3$. Considering \eqref{min_sys_bound} we let $K_0 = 2\times (.9137)^2 = 1.6697$. Then by Proposition \ref{Teo bound} we have
\be \label{ub-s-1/5}
\sup_{z \in X_3}\sum_{i=1}^{3g-3+n}|\mu_{i}|^2(z)\leq \frac{K_0}{\lsys}=\frac{1.6697}{\lsys}.
\ene

On the thick part of the moduli space $\sM^n_g$, the \wep curvature has been well studied in \cite{Huang, Teo09, Wolf-W-1}. Now we study the \wep curvatures on Riemann surfaces with short systoles. Our first result in this section is as follows.
\bt[=Theorem \ref{i-LUB-Ric}] \label{LUB-Ric}
For any $X \in \sM_g^n$ with $\lsys \leq 2\eo$, then
\ben
\item for any $\mu \in T_X\sM_g^n$ with $||\mu||_{WP}=1$, the Ricci curvature satisfies 
\[\Ric(\mu)\geq -\frac{4}{\lsys}.\] 

\item  The scalar curvature at $X$ satisfies 
\[\Sca(X) \geq -\frac{4}{\lsys}\cdot (3g-3+n).\]
\een 
\et
\bp
We first show Part (1). Let $\mu \in T_X\sM_g^n$ with $||\mu||_{WP}=1$ and one may choose a holomorphic orthonormal basis $\{\mu_{i}\}_{i=1}^{3g-3+n}$ of $T_X\sM_g$ such that $\mu=\mu_1$. Now we split the lower bound in \eqref{RicK-eq} into three parts. Since $X_1, X_2$ and $X_3$ are mutually disjoint,
\beqar
\Ric(\mu) &\geq& -2 \sum_{j=1}^{3g-3+n}\int_{X} D(|\mu|^2)\cdot (|\mu_{j}|^2)dA \nonumber\\
&=& -2\int_{X_1} D(|\mu|^2)\cdot ( \sum_{j=1}^{3g-3+n}|\mu_{j}|^2)dA \\
&& -2\int_{X_2} D(|\mu|^2)\cdot ( \sum_{j=1}^{3g-3+n}|\mu_{j}|^2)dA \\
&& -2\int_{X_3} D(|\mu|^2)\cdot ( \sum_{j=1}^{3g-3+n}|\mu_{j}|^2)dA. 
\eeqar
Since $D$ is a positive operator (see \cite{Wol86}), $D(|\mu|^2)\geq 0$. Then it follows by \eqref{ub-l-1/5}, \eqref{ub-l-1/5-cusp} and \eqref{ub-s-1/5} that
\beqar
\Ric(\mu) &\geq& -2\cdot C(\eo)^2\cdot \int_{X_1\bigcup X_2} D(|\mu|^2)dA  -2\cdot \frac{K_0}{\lsys}\cdot \int_{X_3} D(|\mu|^2)dA \\
&\geq & -\frac{3.3394}{\lsys} \int_X D(|\mu|^2)dA
\eeqar
where in the last inequality we note that $2K_0 = 4(.9137^2) = 3.3394$ and $C(\eo) = .7438$.
 Recall that the operator $D$ is self-adjoint and $D(1)=1$. So $\int_X D(|\mu|^2)dA=\int_X |\mu|^2\cdot D(1)dA=||\mu||_{WP}^2=1$. Therefore 
\bear \label{0-ric-lb}
\Ric(\mu) \geq-\frac{3.3394}{\lsys}\geq -\frac{4}{\lsys}.
\eear

Part (2)  follows by Part (1) as
\bear
\Sca(X) =\sum_{i=1}^{3g-3+n}\Ric (\mu_i) \geq -\frac{4}{\lsys}\cdot (3g-3+n).
\eear

The proof is complete.\ep

\br
For $\sM_g=\sM_g^0$, the lower bound in Part (2) of Theorem \ref{LUB-Ric} can be extended to $-\frac{11}{\lsys}\cdot (g-1)$ because \eqref{0-ric-lb} implies that
\beqar
 \Sca(X) &=&\sum_{i=1}^{3g-3}\Ric (\mu_i) \geq -\frac{-3 \times 3.3394}{\lsys}\cdot (g-1)\\
&\geq& -\frac{11}{\lsys}\cdot (g-1).
\eeqar
\er

Since the \wep sectional curvature is negative \cite{Tromba86, Wol86}, we have that for any $X\in \sM_g^n$ and $\mu, v \in T_X\sM_g^n$,
\begin{equation}\label{AllK-eq}
\max\{ \Ric(\mu), \Ric(v)\}<K^{WP}(\mu,v).
\end{equation}
The following result is a direct consequence of Theorem \ref{LUB-K}.
\bt \label{LUB-K}
For any $X \in \sM_g^n$ with $\lsys\leq2\eo$, then for any $\mu,v \in T_X\sM_g^n$, the \wep sectional curvature satisfies that
\bear \label{ulb-K}
K^{WP}(\mu,v)\geq -\frac{4}{\lsys}.
\eear
\et

\br
Huang in \cite{Huang-ajm} showed that $K^{WP}(\mu,v)\geq -\frac{c(g)}{\lsys}$ on $\sM_g$ where $c(g)>0$ is a constant depending on $g$.
\er

\br
The upper bound $2\eo$ for $\lsys$ in Theorem \ref{LUB-K} may not be optimal. However, the upper bound for $\lsys$ can not be removed: actually it was shown in \cite[Theorem 1.1]{Wolf-W-1} that if $\lsys$ is large enough, then $$\min_{\emph{span}\{\mu,v\} \subset T_X\sM_g} K^{WP}(\mu,v)\leq -C<0$$ where $C>0$ is a uniform constant independent of $g$. In particular, \eqref{ulb-K} does not hold for Buser-Sarnak surface $\mathcal{X}_g$ in \cite{BS94} whose injectivity radius grows like $\ln{(g)}$ as $g\to \infty$. 
\er

We close this subsection by proving Theorem \ref{i-ub-orth}.
\bt [=Theorem \ref{i-ub-orth}] \label{ub-orth}
For any $X \in \sM_g^n$ with $\lsys \leq 2\eo$, then for any $\mu\neq 0 \in P(X)^\perp$ and $v\in T_X \sM_g^n$, the \wep sectional curvature $K^{WP}(\mu,v)$ along  then plane spanned by $\mu$ and $v$ satisfies that
\beqar \label{ulb-K}
K^{WP}(\mu,v)\geq -4.
\eeqar
\et

\bp
Since $\mu \in P(X)^\perp$, by Lemma \ref{ub-v-orth} we have \[\sup_{z\in X}|\mu|(z)\leq \sqrt{2}||\mu||_{WP}.\] 
By taking a rescaling one may assume $||\mu||_{WP}=1$. We normalize $v$ such that $||v||_{WP}=1$. Then it follows by \eqref{SecK-eq} that
\bear
K^{WP}(\mu,v) &\geq&  -2\int_X D(|v|^2)|\mu|^2dA \\
&\geq & -4\int_X D(|v|^2)\cdot 1 dA \nonumber \\
&=&-4\int_X |v|^2 dA=-4 \nonumber
\eear
which completes the proof.
\ep

\section{Total scalar curvature for large genus}\label{s-WP-scal}
It is known \cite{Schumacher, Wolpert5} that the \wep scalar curvature always tends to negative infinity as the surface goes to the boundary of the moduli space. In this section we focus on $\sM_g$ and study the total \wep scalar curvature $\int_{\sM_g} \Sca(X)dX$ over the moduli space $\sM_g$, where $dX$ is the \wep measure induced by the \wep metric on $\sM_g$.

For any $\epsilon>0$, the $\epsilon$-thick part $\sM^{\geq \epsilon}_g$ is the subset defined as $$\sM^{\geq \epsilon}_g:=\{X\in \sM_g: \ \lsys\geq \epsilon\}.$$ The complement $\sM^{< \epsilon}_g:=\sM_g \setminus \sM^{\geq \epsilon}_g$ is called the $\epsilon$-thin part of the moduli space. We first recall the following result of Mirzakhani which we will apply.
\bt (Mirzakhani, \cite[Corollary 4.3]{MM4}) \label{mm-1/sys}
As $g\to \infty$,
\[\int_{\sM_g}\frac{1}{\lsys}dX \asymp \Vol(\sM_g).\]
\et

Now we are ready to state our result in this section.
\bt[=Theorem \ref{i-total}] \label{total}
As $g\to \infty$,
\[\frac{\int_{\sM_g} \Sca(X)dX}{\Vol(\sM_g)} \asymp -g.\]
\et

\bp
First by Wolpert \cite{Wol86} or Tromba \cite{Tromba86} we know that for all $X\in \sM_g$,
\[\Sca(X)\leq \frac{-3}{4\pi}\cdot (3g-2).\]
Thus, 
\bear\label{ub-tsca}
\frac{\int_{\sM_g} \Sca(X)dX}{\Vol(\sM_g)} \leq -C_1 \cdot g
\eear
where $C_1>0$ is a uniform constant independent of $g$.

Next we prove the other direction. That is to show that
\bear \label{ulb-tsca}
\int_{\sM_g} \Sca(X)dX \geq -C_1'\cdot g \cdot \Vol(\sM_g)
\eear
where $C_1'>0$ is a uniform constant independent of $g$. We split the total scalar curvature into two parts. More precisely we let $\eo=\sinh^{-1}(1)>0$,
\bear \label{eq-1/5}
\quad \ \ \ \ \int_{\sM_g} \Sca(X)dX=\int_{\sM^{\geq \eo}_g} \Sca(X)dX+\int_{\sM^{<\eo}_g} \Sca(X)dX.
\eear

 On $\sM^{\geq \eo}_g$ it follows by Lemma \ref{teo} of Teo that
\be
\Sca(X)\geq -(6g-6) \cdot C^2(\eo). \nonumber
\ene
Thus, we have
\bear \label{1/5-1}
\int_{\sM^{\geq \eo}_g} \Sca(X)dX &\geq& -(6g-6)\cdot C^2(\eo)\cdot \Vol(\sM^{\geq \eo}_g)  \\
&\geq &-(6g-6)\cdot C^2(\eo) \cdot \Vol(\sM_g) \nonumber\\
&\geq & -C_2\cdot g \cdot \Vol(\sM_g)\nonumber
\eear
where $C_2>0$ is a uniform constant independent of $g$.

On $\sM^{< \eo}_g$ it follows by Theorem \ref{LUB-Ric} that
\be
\Sca(X)\geq -\frac{11}{\lsys}\cdot (g-1). \nonumber
\ene
Thus, we have
\bear 
\int_{\sM^{< \eo}_g} \Sca(X)dX &\geq& -11(g-1)\cdot \int_{\sM^{< \eo}_g}\frac{1}{\lsys}dX \nonumber \\
&\geq & -11(g-1)\cdot \int_{\sM_g}\frac{1}{\lsys}dX. \nonumber 
\eear
By Theorem \ref{mm-1/sys} of Mirzakhani we have
\be \label{1/5-3}
\int_{\sM^{< \eo}_g} \Sca(X)dX \geq -C_3\cdot g \cdot \Vol(\sM_g)
\ene
where $C_3>0$ is a uniform constant independent of $g$.

Then the claim \eqref{ulb-tsca} follows by \eqref{eq-1/5}, \eqref{1/5-1} and \eqref{1/5-3}.
\ep

\end{document}